%% file: main.tex
\documentclass[preprint,12pt]{elsarticle}

\usepackage{amssymb}
\usepackage{amsmath}
\usepackage{amsthm}

\theoremstyle{plain}
\newtheorem{thm}{Theorem}
\theoremstyle{remark}
\newtheorem{rem}{Remark}
\theoremstyle{definition}
\newtheorem{dfn}{Definition}

\usepackage{lineno}

\usepackage{bm}
\usepackage{color}
\usepackage{graphicx}
\usepackage{algorithm}
\usepackage{algorithmicx}
\usepackage{algpseudocode}

\include{abbr-notation}

\newcommand{\bff}[0]{\vec{f}}
\renewcommand{\bu}[0]{\vec{u}}
\renewcommand{\bv}[0]{\vec{v}}
\renewcommand{\bw}[0]{\vec{w}}

\newcommand{\obf}[0]{\overline{\vec{f}}}
\newcommand{\obu}[0]{\overline{\vec{u}}}
\newcommand{\obv}[0]{\overline{\vec{v}}}
\newcommand{\obw}[0]{\overline{\vec{w}}}

\newcommand{\laxparL}[0]{\varepsilon}
\newcommand{\laxparR}[0]{\theta}

\journal{Journal of Computational and Applied Mathematics }

\renewcommand{\vec}[1]{\bm{#1}}

\renewcommand{\R}[0]{\mathbb{R}}

\newcommand{\revAS}[1]{#1} 

\newcommand{\revASi}[1]{#1} 

\begin{document}

\begin{frontmatter}



\title{Coupled general Riemann problems for the Euler equations} 


\author[aff-iapcm]{Zhifang Du}
\affiliation[aff-iapcm]{organization={Institute of Applied Physics and Computational Mathematics},
  postcode={100094}, 
  city={Beijing},
  country={China}}
\ead{du_zhifang@iapcm.ac.cn}

\author[aff-rwth]{Aleksey Sikstel\corref{cor1}}
\affiliation[aff-rwth]{organization={RWTH Aachen University},
  addressline={Schinkelstrasse 2},
  postcode={52062},
  city={Aachen},
  country={Germany}}
\ead{sikstel@acom.rwth-aachen.de}
\cortext[cor1]{Corresponding Author}

\begin{abstract}
We introduce a novel method for systems of conservation laws coupled at a sharp interface based on \revAS{generalized} Riemann problems. This method yields a piecewise-linear in time approximation of the solution at the interface, thus, descynchronising the solvers for the coupled systems. We apply this framework to a problem of compressible Euler equations coupled via a gas generator and prove its solvability. Finally, we conduct numerical experiments and show that our algorithm performs at correct convergence rates.
\end{abstract}


\begin{keyword}
general Riemann problem \sep system of conservation laws   \sep coupling conditions \sep sharp interface \sep compressible Euler equations  \sep gas networks 

35L60 \sep 65M60 \sep 74J40 \sep 90B10  
\end{keyword}

\end{frontmatter}



\input{intro}
\input{coupling-GRP}

\input{grp_analysis}
\input{solvability}

\input{num_results}



\bibliographystyle{elsarticle-num} 
\bibliography{extracted.bib} 


\end{document}

%% file: abbr-notation.tex
\def\bu{\mathbf{u}}
\def\bv{\mathbf{v}}
\def\bw{\mathbf{w}}

\def\pt{\partial}

\def\f#1#2{\frac {#1}{#2}}
\def\f32{\frac 32}
\def\R{{\bf R}}

\def\beq{\begin{equation}}
\def\eeq{\end{equation}}
\def\bga{\begin{array}}
\def\eda{\end{array}}

\def\dfr#1#2{\displaystyle{\frac{#1}{#2}}}

%% file: intro.tex
\section{Introduction}

\begin{figure}[!htb]
  \centering
  
  \includegraphics[width=0.6\textwidth]{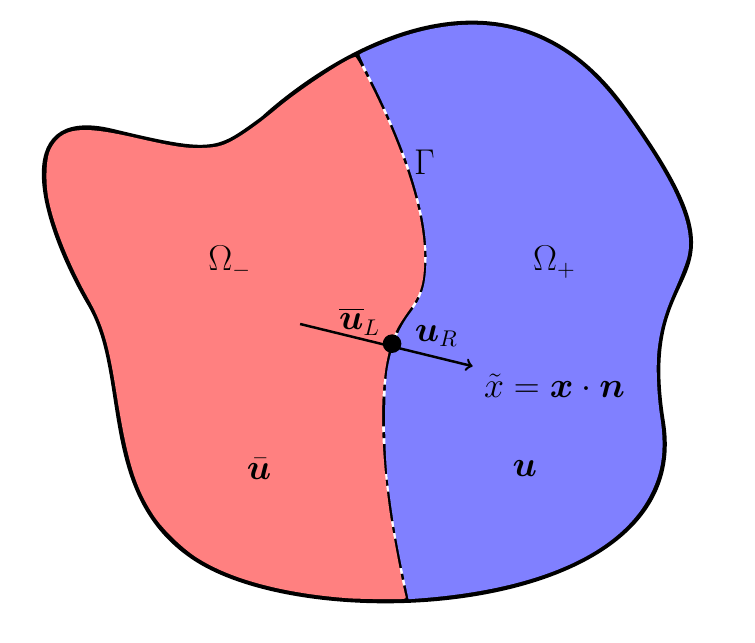}
  \caption{General setup of two coupled hyperbolic conservation laws.}
  \label{fig:general-cpl-setup}
\end{figure}

Coupled systems of hyperbolic conservation laws arise from questions related
to modelling of fluid-structure interaction, pipeline networks, traffic flow, as well as bio-medical, manufacturing and logistic
systems, see~\cite{bressanFlowsNetworksRecent2014} for an overview. In this work, we consider the dynamics
of two systems of conservation laws 
\begin{align}
  \label{eq:pde-left}
  \begin{split}
    &\obu_t + \mathrm{div}\left( \obf (\obu) \right) = \vec{0}, \quad \obu \,\colon\, [0, T] \times \Omega_- \to \overline{\mathcal{D}}\subset \R^{\overline{m}}\\
    &\obu(0, \vec{x}) = \obu_0(\vec{x}),\quad \vec{x} \in \Omega_-
  \end{split},
\end{align}
and
\begin{align}
  \label{eq:pde-right}
  \begin{split}
    &\bu_t + \mathrm{div}\left( \bff (\bu) \right) = \vec{0}, \quad  \bu\,\colon\,[0, T] \times \Omega_+ \to \mathcal{D}\subset  \R^m\\
    &\bu(0, \vec{x}) = \bu_0(\vec{x}),\quad \vec{x} \in \Omega_+
  \end{split},
\end{align}
defined on subdomains $\Omega_- \cup \Omega_+ \cup \Gamma = \Omega \subset \R^d $ separated by an interface $\Gamma$, as depicted in Figure~\ref{fig:general-cpl-setup}. The flux functions $\overline{\vec{f}}$ and $\vec{f}$ map elements of the respective admissible domains $\overline{\mathcal{D}}$ and $\mathcal{D}$ to $\R^{\overline{m}}$ and $\R^m$.
On the interface $\Gamma$ coupling conditions are prescribed in algebraic form, i.e.  $\Psi(\obu_L, \bu_R) =\vec{0}$ is required to hold at the interface for almost all times $t\in [0,T]$ in  normal direction of $\Gamma$, \revAS{where $\Psi \,\colon\, \R^{\overline{m}} \times \R^m \to \R^{\max\{\overline{m}, m\}} $}.  In \revAS{order to distinguish} quantities defined on $\Omega_+$, we denote quantities defined on $\Omega_-$  with an over-line~$\overline{\,\cdot\,}$.

First approaches to the numerical solution of the coupled problem have been iterative  where at each timestep
the conservation laws are solved alternately as in \cite{dickoppCouplingElasticSolids2013, houNumericalMethodsFluidStructure2012}. Such methods are computationally expensive and, typically, lack proofs of convergence  of the
iterative method. In order to avoid these problems,  solutions of coupled half-Riemann problems (half-RPs) have been developed e.g.~in~\cite{herty2008coupling, brouwerGasPipelineModels2011, colomboCauchyProblemPsystem2008, gugatCouplingConditionsTransition2017, herty2021coupling}, for a more comprehensive review on the literature  we refer to~\cite{sikstelAnalysisNumericalMethods2020}. The idea of these methods is to solve simultaneously two RPs, i.e. Cauchy problems for piecewise constant initial conditions with a single jump,  such that the middle state fulfils the coupling conditions. The jump of each half-RP is positioned at $\Gamma$ and the (coupled) intermediate values are used to compute fluxes in the spirit of Godunov-type methods. 

The half-RP methods yield boundary conditions for both systems of conservation laws at the interface $\Gamma$ at each discretisation point in space-time. Thus, no iteration is required and two separate solvers with explicit time-discretisation schemes can be employed  directly. Nonetheless, the time-discretisation schemes for the two systems need to march synchronously. In case the characteristic speeds, i.e.~the maximal stable timestep size, of the two systems differ by a large value, \revAS{synchronisation causes an} unnecessary fine timestep in the slower flow. This drawback can be treated by multilevel time stepping,~\cite{hertyFluidStructureCouplingLinear2018a}, that, however, requires rather complicated bookkeeping and is cumbersome to implement in parallelised solvers.

In this work we present a novel approach using half \emph{generalized} Riemann problems (half-GRPs). A one-dimensional generalized Riemann problem for a conservation law is a  Cauchy problem where the initial condition is a piecewise \emph{polynomial} function with a single jump \cite{BF1984,BF2003,ben2006direct,BL2007}. Since the solution, besides the jump, is not constant in space, the characteristic curves are not necessary linear and the intermediate value of the resulting solution $\vec{w}$, i.e.~$\vec{w}(t, \vec{x} \in \Gamma)$, might change. Besides other quantities, one obtains the partial time derivative in the point $(t=0, \vec{x} \in \Gamma)$ that yields a linear function for the coupled intermediate values, $\vec{w}(t, \vec{x} \in \Gamma)$.

Our idea consists in solving two \revASi{(nonlinear)} half-GRPs such that at every point in $\Gamma$ their solutions fulfil the coupling conditions  $\Psi(\obu_L(t), \bu_R(t)) = \vec{0}$ for any $t$ in a time-interval. Since the coupling conditions are prescribed in the direction normal to $\Gamma$ the coupling problem is inherently one-dimensional. Thus, the analysis of the time-derivatives can rely on the rich amount of previous results, e.g.~\cite{ben2006direct,BL2007,axis-symmetry-grp,DL2020,thermo-grp-2017}. The second order approximation of the solution in time at the interface allows to desynchronise the solvers and, as long as the characteristics curves do not cross the time axis a second time, fully decouples the two solvers. Thus, the half-GRP framework allows more flexibility (e.g. different time discretisation schemes) and facilitates the implementation of large-scale parallel codes on distributed memory systems without tedious multilevel time stepping.

\revASi{A related approach based on higher-order ADER schemes has been proposed in~\cite{borsche2014ader, muller2015high, contarino2016junction}. In contrast to our ansatz, ADER schemes employ linearized Lax-curves at the interface to solve the coupling problem. Thus, the higher-order derivatives at the coupling interface can always be obtained as long as the coupling conditions are sufficiently smooth. However, despite decent numerical results, the linearized setting is not suitable to investigate the  existence of a unique high-order solution to coupled \emph{nonlinear} conservation laws analytically.  }

The outline of this text is as follows. In Section 2 we introduce the framework of coupled GRPs and summarise our findings in an abstract algorithm. In the subsequent Section we apply the GRP framework to the gas-generator coupling problem considered in~\cite{hertyCouplingCompressibleEuler2019} \revAS{ and analyse its solvability}. Finally, in the fourth Section we present numerical results,  measure empirical order of convergence and conclude with a discussion of the results.


%% file: coupling-GRP.tex
\section{Coupled generalized half Riemann problems }
\label{sec:coupl-hyperb-cons}
In the following we develop a framework of coupled GRP methods. We recall previous results on coupled RPs, i.e.~solutions to the coupling problem that are constant in time. These results are extended by including time derivatives of the coupled solutions at the interface that we obtain from half--GRP solutions.

\begin{dfn}
  \label{def:cpl-condition} 
  Let a domain $\Omega \subset \R^d$  be separated by a smooth interface $\Gamma$, i.e. $\Omega =
  \Omega_- \cup \Omega_+ \cup \Gamma$, as illustrated in Figure \ref{fig:general-cpl-setup}.  Furthermore, let $(\obu(t, \vec{x}),\bu(t, \vec{x}) )$ be a pair of weak entropy solutions  to the system of equations~\eqref{eq:pde-left} and~\eqref{eq:pde-right}, respectively. The coupling conditions are prescribed by a sufficiently smooth function $\revAS{\Psi \, \colon\, \R^{\overline{m}}\times \R^m \to \R^{\max\{\overline{m}, m\}}}$.

  Assuming the traces of $\obu$  and $\bu$ at $\Gamma$ exist, we call $(\obu(t, \vec{x}),\bu(t, \vec{x}) )$ a \emph{weak solution of the coupled problem} if  a.e.~in $t$ the coupling condition $\Psi = \vec{0}$ holds. To be more precise, the solution pairs satisfy
  \begin{equation} 
    \label{eq:cpl-conditions}
    \Psi ( \obu (t, \vec{x}^-), \bu(t, \vec{x}^+)) = \vec{0}, \, \text{ for }t > 0, \vec{x}\in\Gamma,
  \end{equation}
  in the limits $\cdot (t, \vec{x}^{\pm}) = \lim_{\varepsilon \downarrow 0} \cdot (t, \vec{x} \pm \varepsilon\vec{n})$ and normal direction $\vec{n} = \vec{n}(\vec{x})$ of the interface $\Gamma$ at $\vec{x}$.
\end{dfn}
The first step of the construction of  coupled weak solutions are coupled half--RPs, investigated e.g.~in~\cite{colomboCauchyProblemPsystem2008, hertyFluidstructureCouplingLinear2018}. We project the conservation laws~\eqref{eq:pde-left} and~\eqref{eq:pde-right} to the normal direction of $\Gamma$ and define two half--RPs that are coupled by conditions $\Psi$ at $\Gamma$. A half--RP in direction $\vec{n} \in\R^d \backslash \{\vec{0}\}$ for the conservation law~\eqref{eq:pde-left} on $\Omega_-$  is given by
\begin{equation}\label{eq:half-rp-left}
  \obu_t + (\overline{\vec{f}}(\obu, \vec{n}))_{\tilde{x}} = \vec{0},
  \quad
  \obu(0, \tilde{x}) = 
  \begin{cases}
    \obu_L &\text{ if } \tilde{x} < 0,\\
    \obu_{\Gamma} &\text{ if } \tilde{x} > 0
  \end{cases},
\end{equation}
where $\obu_L$ and $\obu_{\Gamma}$ are constant and  $\tilde{x}:= \vec{x} \cdot \vec{n}$. The solution $\obu$ of the half--RP is the solution  of~\eqref{eq:half-rp-left}  {\em restricted to} $\Omega_-$. Similarly, we consider
a half--RP associated with the projected conservation law~\eqref{eq:pde-right} on $\Omega_+$, i.e.~we solve
\begin{equation}\label{eq:half-rp-right}
  \bu_t + (\vec{f}(\bu, \vec{n}))_{\tilde{x}} = \vec{0},
  \quad
  \bu(0, \tilde{x}) = 
  \begin{cases}
    \bu_{\Gamma} &\text{ if } \tilde{x}  < 0,\\
    \bu_R &\text{ if } \tilde{x}  > 0
  \end{cases},
\end{equation}
 where again $\bu_{\Gamma}$ and $\bu_R$ are constant. Particularly interesting  are  solutions  \revAS{whose trace}  in $\Omega_-$,  $\obu(t,\vec{x}^-)=\overline{\bv}$, and the trace in $\Omega_+$, $\vec{w}(t,\vec{x}^+)=\bv$, fulfil the coupling conditions $\Psi(\obv, \bv)=\vec{0}$.

Given a constant state $\obu_L$ we  introduce  admissible boundary
states  as follows. The set
\begin{align}\label{eq:adm-set-half-rp1}
  \begin{split}
  \overline{V}(\obu_L) := \bigg\{ \overline{\vec{v}} \in \overline{\mathcal{D}} \, \colon  \, &\exists \overline{M}\in\{1,\ldots,\overline{m}\},\, \exists\, (\laxparL_i)_{i=1}^{\overline{M}}\subset\R^{\overline{m}} \text{ s.t. } \obv = \obv_{\overline{M}} \text{ with } \\
  & \obv_0 := \obu_L,\, \obv_i := \overline{L}_i^+(\laxparL_i; \obv_{i-1}),\,  i=1,\ldots,\overline{M} \text{ and}\\
    & \overline{\lambda}_{\overline{M}}(\obv)  \leq 0 < \overline{\lambda}_{\overline{M}+1}(\obv)   \bigg\}
   \end{split}
\end{align}
consists of all states $\overline{\bv} = \obu(t,\vec{x}^-)$ that solve a  half-RP in $\Omega_-$, i.e.~are attainable by a composition of forward Lax curves $\overline{L}^+_{i}(\cdot \, ;\, \obu^{ \, i-1 }_0),\, i \in \{1,\ldots\overline{M}\}$ for some $\overline{M}\in\{1,\ldots,\overline{m}\}$ each emerging at some $\obu^{i-1}_0 \in \overline{\mathcal{D}}$ located on the previous Lax curve. In particular, $\overline{M}$ is chosen such that for the eigenvalues of the Jacobian of the flux $\overline{\vec{f}}$ in~\eqref{eq:pde-left} $\overline{\lambda}_{\overline{M}}(\obv)\leq 0$ and $\overline{\lambda}_{\overline{M}+1}(\obv)>0$ holds. Here we use the convention $\lambda_{\overline{m} + 1}  \equiv +\infty$. The composition of these Lax curves connects $\vec{w}_0$ to a state that is located at the interface.
Similarly, for a given state
$\vec{w}_R$ the set of admissible boundary states in $\Omega_+$ is defined by
\begin{align}\label{eq:adm-set-half-rp2}
  \begin{split}
  V(\bu_R) := \bigg\{ \vec{v} \in \mathcal{D} \, \colon \, &\exists M\in\{1,\ldots,m\},\, \exists\, (\laxparR_i)_{i=1}^{M}\subset\R^{m} \text{ s.t. } \bv = \bv_M \text{ with } \\
  &\bv_0 := \bu_R,\, \bv_i := L_i^-(\laxparR_i; \bv_{i-1}),\,  i=1,\ldots,M \text{ and}\\
   & \lambda_{M-1}(\bv)  \leq 0 < \lambda_{M}(\bv)   \bigg\}
 \end{split}
\end{align}
where $M$ is chosen such that the eigenvalues of the Jacobian of the flux $\vec{f}$ in~\eqref{eq:pde-right} $\lambda_M(\bv) > 0$ and $\lambda_{M-1}(\bv) \leq \, 0$. Here we use the convention $\lambda_0 \equiv -\infty$.
\begin{dfn}
  \label{def:coupl-riem-probl}
Provided there exist unique states $\overline{\bv}\in \overline{V}(\obu_L)$ and $\bv \in V(\vec{u}_R)$
such that $\vec{\Psi}( \overline{\bv}, \bv ) = \vec{0}$, the \emph{solution to the coupled RP}~\eqref{eq:half-rp-left},~\eqref{eq:half-rp-right} is given by the solution of the 
two half--RPs with respective initial data $(\obu_L, \obu_{\Gamma} )$ 
and $( \bu_{\Gamma}, \bu_R)$.
\end{dfn}
The solution of the coupled RP can  be employed to solve the coupled problem, e.g.~by setting boundary conditions for each system at time points such as Runge-Kutta stages. This approach allows to use separate solvers for each PDE, however requires full synchronisation of the time discretisation. 

Instead of enforcing coupling conditions at discrete points in time by means of half--RPs, we introduce continuous-in-time,  second order  approximations of the values coupled at the interface $\Gamma$. To this end, we consider a time interval $[t_0, t_1]$, denote the initial data at $\Gamma$ by $\obu_L := \obu(t_0, \vec{x}^-)$ and $\bu_R := \bu(t_0, \vec{x}^+)$. Moreover, we henceforth call the coupled RP the \emph{associated coupled RP} and denote its solutions by $(\cdot)_*$ and $\overline{(\cdot)}_*$. Our aim is the solution of the coupled GRP, i.e.~we search for two half--GRP solutions that fulfil the coupling conditions:
\begin{align}
  \label{eq:grp-sol-def}
  \begin{split}
    &\obu_*(t, \vec{x}^-) = \obu_* + (t-t_0)\left(\frac{\partial \obu}{\partial t}\right)_* , \quad
      \bu_*(t, \vec{x}^+) = \bu_* + (t-t_0)\left(\frac{\partial \bu}{\partial t}\right)_*, \\
    &\text{such that } \Psi(\obu_*(t, \vec{x}^-), \bu_*(t, \vec{x}^+))  = \vec{0} \text{ for } t \in [t_0, t_1],
  \end{split}
\end{align}
where the derivatives marked with $(\cdot)_*$ are defined as 
\begin{equation*}
  \left(\frac{\partial \obu}{\partial t}\right)_* := \lim_{t\downarrow t_0}\frac{\partial \obu}{\partial t}(t, \vec{x}^-) \ \ \text{ and } \ \ \left(\frac{\partial \bu}{\partial t}\right)_* := \lim_{t \downarrow t_0}\frac{\partial \bu}{\partial t}(t, \vec{x}^+).
\end{equation*}
In addition, we require the time interval $[t_0, t_1]$ to be small enough such that no wave crosses the $t$-axis \revAS{two times}.

Next, we characterise the sets of the time derivatives  compatible with the solutions of the associated coupled half-RPs. The solution of a single GRP involves a transformation of the PDE onto characteristic coordinates and analysis of the Riemann invariants, cf.~\cite{ben2006direct,thermo-grp-2017}. Typically, this procedure boils down to  a linear algebraic system where the coefficients and the right hand side depend on the initial states and the middle state  of the associated RP: 
\begin{equation}
  \label{eq:grp-sol-general}
  \vec{C}(\bu_L, \bu_R, \bu_*)\left( \frac{\partial \bu}{\partial t} \right)_* = \vec{g}(\bu_L, \bu_R,  \bu_*)
\end{equation}
with $\vec{C} \in \R^{m\times m}$ and $\vec{g} \in \R^m$. 

We assume that there exists a unique solution of the associated coupled RP.  Since, each associated RP is, in fact, a half--RP and the time derivatives of the GRP require information of the whole middle state, i.e.~for $\vec{x} = \vec{x}^{\pm}$, the coefficients of the coupled GRP depend on both $\obu_*$ and $\bu_*$.  Thus, the set of admissible time derivatives at the interface $\Gamma$ read
\begin{align}
  \label{eq:adm-derivatives-set-left}
  \begin{split}
    \overline{V}'(\obu_L,  \obu_*, &\bu_*) \\
    := \big\{ &\obv' \in \R^{\overline{m}} \,\colon\,\,  \obv'\! =\! \overline{\vec{C}}^{-1}(\obu_L, \bw, \obu_*, \bu_*)\overline{\vec{g}}(\obu_L, \bw, \obu_*, \bu_*), \\
               &\bw \in \mathcal{D},\\
                                         &(\obu_*, \bu_*) \text{ -- middle state of the assoc.~coupled RP},\\                                          &\overline{\vec{C}}\in \R^{\overline{m}\times\overline{m}} \text{ -- matrix of the coefficients for the GRP solution},\\
    &\overline{\vec{g}} \in \R^{\overline{m}} \text{ -- rhs of the GRP solution} \big\},
  \end{split}
\end{align}
and, analogously, for the right-hand side we define
\begin{align}
  \label{eq:adm-derivatives-set-right}
  \begin{split}
    V'(\bu_R,  \obu_*, &\bu_*) \\
    := \big\{ &\bv' \in \R^m \,\colon\,\,  \bv'\! =\! \vec{C}^{-1}(\obw, \bu_R,  \obu_*, \bu_*)\vec{g}(\obw, \bu_R, \obu_*, \bu_*), \\
               &\obw \in \overline{\mathcal{D}},\\
                                         &(\obu_*, \bu_*) \text{ -- middle state of the assoc.~coupled RP},\\                                          &\vec{C} \in\R^{m\times m} \text{ -- matrix of the coefficients for the GRP solution},\\
    &\vec{g} \in\R^m \text{ -- rhs of the GRP solution} \big\}.
  \end{split}
\end{align}
The matrices $\overline{\vec{C}}$ and $\vec{C}$ as well as corresponding right-hand sides $\overline{\vec{g}}$ and $\vec{g}$ can be obtained by applying the GRP derivation procedure to each coupled half--RP. In the next Section we are going to present this derivation for the case of compressible  Euler equations in detail. We summarise the above ideas as a general Algorithm~\ref{alg:CGRP}.
\begin{algorithm}[!h]
  \caption{Generalized Riemann solver for coupled problems.} \label{alg:CGRP}

  \begin{algorithmic}[1]
    \Require
    \Statex current time $t_0$, position at interface $\vec{\xi} \in \Gamma$
    \Statex normal of the interface $\vec{n} = \vec{n}(\vec{\xi})$
    \Statex adjacent states $\obu_{L}=\obu(t_0, \vec{\xi}^-)$ and $\bu_{R}=\bu(t_0, \vec{\xi}^+)$ in normal direction
    \Statex spatial derivatives in normal direction $\nabla^{+}_{\vec{x}} \obu \cdot \vec{n} \big|_{\vec{x} = \vec{\xi}} $ and $\nabla^{-}_{\vec{x}} \bu \cdot \vec{n} \big|_{\vec{x} = \vec{\xi}} $
    \State  Project the systems \eqref{eq:pde-left} and \eqref{eq:pde-right} onto the  interface normal  $\vec{n}$
    \State Solve the associated coupled RP consisting of two half--RPs at $t_0$, i.e.~find states $\obu_*$ and $\bu_*$
    \State Determine the sets of possible derivatives $\overline{V}'$ and $V'$
    \State Find $(\overline{v}'_*, v'_*) \in \overline{V}'\times V'$ such that $\Psi(\obu_* + (t-t_0) \overline{v}_*', \bu_* + (t-t_0) v_*') = 0$ 
    \State \Return $(\obu_*, \bu_*)$ and $(\overline{v}'_*, v'_*)$
  \end{algorithmic}
\end{algorithm}

The returned values of the Algorithm~\ref{alg:CGRP} represent a second order approximation~\eqref{eq:grp-sol-def}, i.e.~time-dependent boundary conditions for each coupled PDE. In particular the time-discretisations of the PDE-solvers are not required to be synchronised anymore.

\begin{rem}
  
By the chain rule the time derivative of the coupling conditions can be written as
\begin{align*}
  \frac{\partial}{\partial t}\Psi(\obu (t, \vec{x}^-), \bu(t, \vec{x}^+)) &= \nabla_{\obu} \Psi(\obu (t, \vec{x}^-), \bu(t, \vec{x}^+)) \cdot \frac{\partial\obu}{\partial t} \Big|_{t=t_0}\\
  &+ \nabla_{\bu} \Psi(\obu (t, \vec{x}^-), \bu(t, \vec{x}^+)) \cdot \frac{\partial\bu}{\partial t} \Big|_{t=t_0}.
\end{align*}
\revAS{For coupling conditions $\Psi$ that are linear in the first and in the second argument},  one can search for solutions of \revAS{order higher than two}, if the corresponding high-order single GRP solutions are known. In practice, the coupling conditions are often linear (e.g. enforcing continuity of quantities at the interface), however higher-order GRP solutions are, in general, technically difficult to obtain.

\end{rem}


%% file: grp_analysis.tex
\section{Coupled generalized Riemann problem for the Euler equations}
\label{sec:coupl-gener-riem}
This chapter canvasses the implementation of Algorithm~\ref{alg:CGRP} for the gas-generator coupling (GGC) problem of the Euler equations. The GGC problem was previously solved in~\cite{hertyCouplingCompressibleEuler2019} by means of  half-RPs.

\subsection{Euler equations coupled by a gas-generator }
\label{sec:coupled-grp-gas}

The GGC-problem consists of two one-dimensional, compressible gas flows each governed by the Euler equations
\begin{align}
  \label{eq:euler-1d}
  \begin{split}
    &\rho_t + (\rho u)_x = 0\\
    &(\rho u)_t + (\rho u^2 + p)_x = 0\\
    &(\rho E)_t + (\rho u (E + p / \rho))_x = 0
  \end{split},
\end{align}
closed by the ideal gas equations of state $p = \rho R_{sgc} T = \rho e(\gamma -1)$ where $R_{sgc}$ denotes the specific gas constant and $\gamma\geq 1$ the heat capacity ratio. The total energy per unit volume is  $E = \rho e + \frac12 \rho u^2$ where $e$ is the internal energy per unit mass.

The coupling condition at $x=0$ is set such that the momentum $\rho u(t,x = 0)$  has a negative jump of magnitude $\mathcal{E}(t) \geq 0$  while the pressure  $p(t,x = 0) $ as well as the temperature $T(t,x  =0)$  are  continuous at the interface, i.e.
\begin{equation}
  \label{eq:ggc-cpl-conditions}
    \overline{p} =p, \quad
    \overline{\rho} =\rho, \quad
    \overline{u} = u +  \frac{\mathcal{E} }{\overline{\rho}}.
  \end{equation}
  In addition, the solution is required to be subsonic at the coupling interface:
  \[
    |u(t, x = 0^{\pm})| < c(t, x = 0 ^{\pm}) = \sqrt{\frac{\gamma p(t, x = 0 ^{\pm})}{\rho(t, x = 0 ^{\pm})}}.
  \]
  Lastly, we adopt the convention, that the gas flows in positive direction, i.e. $\bar{u}, u \geq 0$. Further details on the solution of the coupled half-RPs can be found in~\cite{hertyCouplingCompressibleEuler2019, sikstelAnalysisNumericalMethods2020}.

\subsection{Generalized Riemann problem for single Euler equations}
Let us consider the GRP for a single system of Euler equations~\eqref{eq:euler-1d} that has been studied e.g.~in~\cite{ben2006direct,BL2007,axis-symmetry-grp,thermo-grp-2017}. \revAS{The solution of the associated RP for the ideal Euler equations consists of three waves: genuinely nonlinear 1- and 3-waves and linearly degenerate 2-wave, see~\cite{toroRiemannSolversNumerical2009}. This means, that the 2-wave is a contact discontinuity with slope $u$ while the 1- and 3- wave are either rarefaction or shock waves with slopes $u-c$ and $u+c$, respectively. }

We assume that the $t$-axis is located in the intermediate region\revAS{, i.e. between the 1- and 2-wave or between the 2- and 3-wave.} The time derivative $\left(\dfrac{\partial \vec{u}}{\partial t}\right)_*$ is obtained from the solution of  the two linear equations
\beq\label{eq:euler-grp-lag}
\bga{l}
a_L\left(\dfr{Du}{Dt}\right)_* + b_L\left(\dfr{Dp}{Dt}\right)_* = d_L,\\[2mm]
a_R\left(\dfr{Du}{Dt}\right)_* + b_R\left(\dfr{Dp}{Dt}\right)_* = d_R,
\eda
\eeq
where  $a_{L/R}$, $b_{L/R}$, and $d_{L/R}$ are determined by resolving left and right waves, respectively.
Their specific values only involve the associated Riemann solution, the initial values and the initial slopes of the fluid state. We refer to \cite[Eqs. (E.3) and (E.7)]{axis-symmetry-grp} for the detailed procedure \revAS{of deriving} these coefficients.

By transforming the total derivatives to \revAS{partial time derivatives, we transform the system~\eqref{eq:euler-grp-lag} to}
\beq\label{eq:euler-grp-eul}
\bga{l}
h_L\left(\dfr{\partial u}{\partial t}\right)_* + k_L\left(\dfr{\partial p}{\partial t}\right)_* = q_L,\\[2mm]
h_R\left(\dfr{\partial u}{\partial t}\right)_* + k_R\left(\dfr{\partial p}{\partial t}\right)_* = q_R,
\eda
\eeq
\revAS{following~\cite{BL2007}, where it was shown that}
\begin{equation}\label{eq:coeff-grp-eul}
h_J = {a_J-\rho_{*J}u_*b_J}, \ \
k_J = {b_J-\dfr{u_*}{\rho_{*J}{c_{*J}}^2}a_J}, \ \
q_J = \left(1-\dfr{{u_*}^2}{{c_{*J}}^2}\right)d_J,
\end{equation}
\revAS{with} $J=L,R$ and $u_*$ is the velocity in the region between the left-going 1-wave and the right-going 3-wave. \revAS{Quantities with subscript} $(\cdot)_{L*}$ ($(\cdot)_{R*}$, resp.) \revAS{denote values} in the region between the \revAS{left-going 1-wave} (right-going 3-wave, resp.) and the \revAS{2-wave} contact discontinuity.

The time derivative of $\rho$ satisfies the linear equation
\beq\label{eq:euler-grp-rho}
g_\rho\left(\dfr{\pt\rho}{\pt t}\right)_* + g_u\left(\dfr{\pt u}{\pt t}\right)_* + g_p\left(\dfr{\pt p}{\pt t}\right)_*=f,
\eeq
where $f$, $g_\rho$, $g_p$, and $g_u$ are known coefficients determined by resolving the left or the right wave. See \cite[Eqs. (E.9) and (E.10)]{axis-symmetry-grp} for details.

\subsection{Coupled generalized Riemann problem for the gas-generator}
Next, we formulate the coupled GRP introduced in Algorithm~\ref{alg:CGRP} for the GGC problem using the GRP solution of a single Euler system discussed above. Recall the notation introduced in Section~\ref{sec:coupl-hyperb-cons}: quantities accented with $\overline{(\cdot)}$ \revASi{are located} to the left of the 1-wave, while quantities lacking an accent \revASi{are} to the right of the 3-wave.  Following the derivation of the first equation in~\eqref{eq:euler-grp-eul}, one recovers the time-derivatives to the left of the interface:
\beq\label{eq:ggc-grp-left}
h_L\left(\dfr{\pt\overline{u}}{\pt t}\right)_* + k_L\left(\dfr{\pt\overline{p}}{\pt t}\right)_* = q_L.
\eeq
Similarly, to the right of the interface we have
\beq\label{eq:ggc-grp-right}
h_R\left(\dfr{\pt u}{\pt t}\right)_* + k_R\left(\dfr{\pt p}{\pt t}\right)_* = q_R.
\eeq
The coefficients $h_{L/R}$, $k_{L/R}$, and $q_{L/R}$ are obtained by substituting the associated coupled Riemann solutions $\obu_*$ and $\bu_*$ into \eqref{eq:coeff-grp-eul}.
Time derivatives $\frac{\partial\overline{u}}{\partial t}$ and $\frac{\partial\overline{p}}{\partial t}$ are located on the left of the interface, while $\frac{\partial u}{\partial t}$ and $\frac{\partial p}{\partial t}$  on the right of the interface.

Taking time derivatives of the coupling conditions~\eqref{eq:cpl-conditions}, we have
\beq\label{eq:ggc-conditions-derivatives}
\bga{l}
\dfr{\pt\overline{p}}{\pt t}-\dfr{\pt p}{\pt t}=0,\\[2mm]
\dfr{\pt\overline{\rho}}{\pt t}-\dfr{\pt\rho}{\pt t}=0,\\[2mm]
\dfr{\pt\overline{u}}{\pt t}-\dfr{\pt u}{\pt t}=\dfr{\mathcal{E}^\prime}{\overline\rho}-\dfr{\mathcal{E}}{\overline\rho^2}\dfr{\pt\overline\rho}{\pt t},
\eda
\eeq
which holds along the sharp interface for $t\in[t_0,t_1]$.
Sending ~\eqref{eq:ggc-conditions-derivatives} to $t=t_0$ and substituting~\eqref{eq:euler-grp-rho} into it lead to
\beq\label{eq:gcc-grp-u}
\left(1-\dfr{\mathcal{E}}{\overline{\rho}_*^2}\dfr{\bar{g_u}}{\bar{g_\rho}}\right)\left(\dfr{\pt\overline{u}}{\pt t}\right)_*
-\left(\dfr{\pt u}{\pt t}\right)_*
-\dfr{\mathcal{E}}{\overline{\rho}_*^2}\dfr{\bar{g_p}}{\bar{g_\rho}}\left(\dfr{\pt\overline{p}}{\pt t}\right)_*
=\dfr{\mathcal{E}^\prime}{\overline\rho_*}-\dfr{\mathcal{E}}{\overline{\rho}_*^2}\dfr{\bar{f}}{\bar{g_\rho}},
\eeq
where $\overline{\rho}_*$ is the density located on the left of the interface obtained by solving the associated coupled RP. The values of the coefficients $\bar{g}_{\rho}$, $\bar{g}_{u}$, $\bar{g}_{p}$, and $\bar{f}$ are evaluated using $\obu_*$.
Combining~\eqref{eq:ggc-grp-left},~\eqref{eq:ggc-grp-right}, and ~\eqref{eq:gcc-grp-u}, the solution of the GGC is \revASi{equivalent to} the solution of \beq\label{eq:gcc-grp}
\bga{l}
h_L\left(\dfr{\pt\overline{u}}{\pt t}\right)_* + k_L\left(\dfr{\pt\overline{p}}{\pt t}\right)_* = q_L,\\[2mm]
h_R\left(\dfr{\pt u}{\pt t}\right)_* + k_R\left(\dfr{\pt\overline{p}}{\pt t}\right)_* = q_R,\\[2mm]
\left(1-\dfr{\mathcal{E}}{\overline{\rho}_*^2}\dfr{\bar{g_u}}{\bar{g_\rho}}\right)\left(\dfr{\pt\bar{u}}{\pt t}\right)_*
-\left(\dfr{\pt u}{\pt t}\right)_*
-\dfr{\mathcal{E}}{\overline{\rho}_*^2}\dfr{\bar{g_p}}{\bar{g_\rho}}\left(\dfr{\pt\overline{p}}{\pt t}\right)_*
=\dfr{\mathcal{E}^\prime}{\overline\rho_*}-\dfr{\mathcal{E}}{\overline{\rho}_*^2}\dfr{\bar{f}}{\bar{g_\rho}},
\eda
\eeq
where the first equation in the coupling condition~\eqref{eq:ggc-conditions-derivatives} is used to eliminate $\frac{\partial{p}}{\partial t}$.
The coefficient matrix of the linear system~\eqref{eq:gcc-grp} for the unknown vector $\left[ (\frac{\partial\overline{u}}{\partial t})_*,(\frac{\partial{u}}{\partial t})_*,(\frac{\partial\overline{p}}{\partial t})_*\right]^\top$  is
\beq\label{eq:coefficient-matrix}
\vec{C}:=\left[
\bga{rrr}
h_L \ \  & 0 \ \ & k_L\\
0 \ \ & h_R \ \ & k_R\\
1-\dfr{\mathcal{E}}{\overline{\rho}_*^2}\dfr{\bar{g_u}}{\bar{g_\rho}} \ \
 & -1  \ \ &
-\dfr{\mathcal{E}}{\overline{\rho}_*^2}\dfr{\bar{g_p}}{\bar{g_\rho}}
\eda
\right].
\eeq\\


%% file: solvability.tex
\revAS{The following Theorem characterizes the solvability of the GGC, i.e. of the  solutions to the linear system~\eqref{eq:gcc-grp}.}
\begin{thm}
The system \eqref{eq:gcc-grp} has a unique solution if the left-going 1-wave of the GGC is a rarefaction wave or a \emph{weak} shock, regardless of the type the right-going 3-wave.
\end{thm}
\begin{proof}
In order to analyze the solvability of the GGC derived in the previous section, we verify that the matrix $\vec{C}$ defined in \eqref{eq:coefficient-matrix} is invertible. We proceed in two steps: firstly checking the signs of elements in the first two rows, secondly showing that $\det(\vec{C})$ is not zero.

According to \cite{axis-symmetry-grp}, in case the 1- and 3-waves are both rarefaction waves, we have
\beq
a_L=1, \ \ \ b_L=\dfr{1}{\overline\rho_*\overline{c}_*}, \ \ \ 
a_R=1, \ \ \ b_R=-\dfr{1}{\rho_*c_*}.
\eeq
Consequently, the signs of the coefficients in the first two rows of $\vec{C}$ are immediately identified as
\beq
\begin{aligned}
&h_L=\dfr{\overline{c}_*-\overline{u}_*}{\overline{c}_*}>0,  \ \
k_L =\dfr{\overline{c}_*-\overline{u}_*}{\overline\rho_*\overline{c}_*^2}>0, \\
&h_R=\dfr{c_*+u_*}{c_*}>0,   \ \
k_R=\dfr{c_*+u_*}{\rho_*c_*^2}<0.
\end{aligned}
\eeq

For a left-going 1-shock,
\beq\label{eq:a-b-shock-left}
a_L=1-\overline\rho_*(\sigma-\overline{u}_*)\Phi_1, \ \ 
b_L=-\dfr{\sigma-\overline{u}_*}{\overline\rho\overline{c}_*^2}+\Phi_1,
\eeq
where $\sigma$ is the shock speed and $\Phi_1>0$ is defined in \cite[Eq. (E.4)]{axis-symmetry-grp}. Thus, the coefficients transform to
\begin{align*}
  h_L&= 1-\overline{\rho}_*(\sigma-\overline{u}_*)\Phi_1-\overline{\rho}_*\overline{u}_*\left(-\dfr{\sigma-\overline{u}_*}{\overline{\rho}_*\overline{c}_*^2}+\Phi_1\right)\\
  &=1-\overline{\rho}_*\sigma\Phi_1
  +\dfr{\overline{u}_*(\sigma-\overline{u}_*)}{\overline{c}_*^2}.
  \\
  k_L&=-\dfr{\sigma-\overline{u}_*}{\overline{\rho}_*\overline{c}_*^2}+\Phi_1-\dfr{\overline{u}_*}{\overline{\rho}_*\overline{c}_*^2}[1-\overline{\rho}_*(\sigma-\overline{u}_*)\Phi_1]\\
  &=-\dfr{\sigma}{\overline{\rho}_*\overline{c}_*^2}+\left[1+\dfr{\overline{u}_*(\sigma-\overline{u}_*)}{\overline{c}_*^2}\right]\Phi_1.
\end{align*}
Since the shock is subsonic with respect to the post shock fluid,
$|\sigma-\overline{u}_*|<|\overline{c}_*|$,
and is traveling left with respect to the fluid, i.e.
$\sigma-\overline{u}_* < 0$, we have $\sigma-\overline{u}_* > -\overline{c}_*$.
By the subsonic assumption,
$|\overline{u}_*| < |\overline{c}_*|$, the relation
$|\overline{u}_*(\sigma-\overline{u}_*)| < \overline{u}_*^2$ holds,
which leads to
\beq
1+\dfr{\overline{u}_*(\sigma-\overline{u}_*)}{\overline{c}_*}>0.
\eeq
Employing $\sigma<0$, we arrive at
\beq
h_L>0, \ \ \ k_L>0.
\eeq
For a right-going \revAS{3-shock} it is easy to verify, by a similar procedure, that
\beq
h_R>0, \ \ \ k_R<0.
\eeq

In conclusion, regardless of the wave types, the signs of the coefficients in the first two rows of $\vec{C}$ are
\beq
h_L>0, \ \ \ k_L>0, \ \ \ h_R>0, \ \ \ k_R<0.
\eeq

In case of a \revAS{left-going}  1-rarefaction wave, according to \cite{axis-symmetry-grp},
\beq
\bar{g_\rho}=\overline{c}_*^2, \ \ \ \bar{g_p}=-1, \ \ \ \bar{g_u}=0.
\eeq
This yields the coefficient matrix
\beq\label{eq:coefficient-matrix-rarefaction}
\vec{C} = \left[
\bga{rrr}
h_L \ \  & 0 \ \ & k_L\\
0 \ \ & h_R \ \ & k_R\\
1 \ \
 & -1  \ \ &
-\dfr{\mathcal{E}}{(\overline\rho_*\overline{c}_*)^2}
\eda
\right],
\eeq
whose determinant is 
\beq
\det(\vec{C})
=
-k_Lh_R+h_Lk_R-\dfr{\mathcal{E}}{(\overline\rho_*\overline{c}_*)^2}h_Lh_R<0.
\eeq

Finally, we consider the case of a left-going 1-shock. 
By continuity, for $\forall\delta$ there exists $\varepsilon>0$ such that
\begin{equation*}
\max(|\bar{g_\rho}-\overline{c}_*^2|,|\bar{g_p}+1|,|\bar{g_u}|)<\delta.
\end{equation*}
when
\begin{equation*}
    |(\overline\rho_*,\overline{u}_*,\overline{p}_*)-(\overline\rho_L,\overline{u}_L,\overline{p}_L)|<\varepsilon
\end{equation*}
According to~\cite{Courant-Friedrichs}, the shock branch of the \revAS{genuinely nonlinear} Lax curve is tangent to the isentropic branch at the second order in the phase space. \revAS{In other words, the Lax curve is two times continuously differentiable at the transition point from a shock to a rarefaction wave.} Thus, $\delta<\hat{C}\varepsilon$ for a fixed constant $\hat{C}>0$.
Since the determinant of $\vec{C}$ continuously depends on all its entries, \revAS{there exists an $\varepsilon$ such that $\det(\vec{C}) < 0$}.

This \revAS{implies that the system has a unique solution for  left-going  rarefaction waves or sufficiently weak shock waves.}
\end{proof}

%% file: num_results.tex
\section{Numerical results}
\label{sec:numerical-results}

We conduct numerical experiments in order to investigate the convergence behaviour of the coupled GRP method for the GGC problem. In all of the following cases we set the material parameters in the EoS of the Euler equations to $\gamma=1.4$ and $R_{sgc} = 277.13333$. Moreover, each boundary condition is set to the corresponding initial condition characteristically.

We compute  numerical solutions $\overline{\vec{u}}_h, \, \vec{u}_h$ by means of third-order SSP-Runge-Kutta discontinuous Galerkin (RKDG) schemes \cite{cockburnRungeKuttaDiscontinuous1998}. The coupling procedure is implemented via Algorithm~\ref{alg:CGRP} and the solvers are implemented in the \textsc{Multiwave}-framework, cf.~\cite{gerhardWaveletFreeApproachMultiresolutionBased2021}. In order to avoid spurious oscillations we employ the classical  Cockburn-Shu limiter from~\cite{cockburnRungeKuttaDiscontinuous1998} on characteristic variables.

The domains $\overline{\Omega}$ and $\Omega$   are tesselated by  $2^{L}  \overline{N_0} $ and $ 2^{L}  N_0 $  cells of uniform width, respectively, where $L$ denotes the level of resolution and  $\overline{N_0}$ as well as  $N_0$ the respective number of cells at  level zero. The timestep-size is set according to the CFL-condition from \cite{chalmersRobustCFLCondition2020} as
\begin{equation*}
  \Delta t = C_{cfl} \frac{\Delta x}{(2P  + 1)\lambda_{max}}
\end{equation*}
with cell width $\Delta x$, maximum characteristic speed $\lambda_{max}$ and polynomial order $P=3$. 

Since there are no smooth analytic solutions to the GGC problem  available, we compute for each test case a finely-resolved numerical solution on level $L=10$. For the discretisation we use again the third-order RKDG scheme except that the states at the coupling boundary are obtained via classical RPs, i.e.~coupled half-RPs are solved at each Runge-Kutta stage as in previous works.  These numerical solutions, denoted by $\overline{\vec{u}}^*, \, \vec{u}^*$, serve as exact solutions  in the computation of the empirical  convergence order (EoC) of the $L^1$-error at final time $T$, i.e.
\begin{equation}
  \label{eq:L1-err}
  err := \left| \|\overline{\vec{u}}_h(T, \cdot) - \overline{\vec{u}}^*(T, \cdot) \|_{L^1(\overline{\Omega})} + \|\vec{u}_h(T, \cdot) - \vec{u}^*(T, \cdot) \|_{L^1(\Omega)} \right|_{\infty}.
\end{equation}

In the first test case we compare solutions obtained by our method with those obtained by the classical RP coupling algorithm. Since the solution features kinks we cannot expect convergence rates higher than 2. In the second one the solution is initially constant and a smooth perturbation in the outtake induces variations in the solution through the  coupling interface. In the third test case a bump in the density of the initial data is transported over the coupling interface that has a constant jump in the momentum.  In both cases we expect a convergence rate of at least two, since the approximation in time of the solution at the coupling interface is of second order. 

\paragraph{Case 1: Linear Outtake}

For the first test case we compare solutions obtained by Algorithm~\ref{alg:CGRP} with the test case of gradually varying, continuous  outtake proposed in~\cite[Sec 4.1]{hertyCouplingCompressibleEuler2019}. The initial density is set to $\rho(0,x)=1$, the velocity to $v(0,x)=1$ and the pressure to $p(0,t)=146820.4$ in both pipes. The outtake function and its derivatve are given by
\begin{equation*}
  \mathcal{E}(t) =
  \begin{cases}
    \min(0.6,\, 3t) &\text{ if } t \leq 3,\\
    \max(0,\, -3t + 1.5) &\text{ otherwise}
  \end{cases}, \quad \mathcal{E}'(t) =
  \begin{cases}
     3   &\text{ if } 0 < t \leq 0.2,\\
    -3  &\text{ if }  0.3 < t \leq 0.5,\\
    0 &\text{ otherwise}
  \end{cases}.
\end{equation*}
The setup of the domain and its discretization parameters are summarized in Table~\ref{tab:params-case-1}. 
\begin{table}[!htb] 
  \centering 
  \begin{tabular}{c|c|cc|cc}
    $T$ & $C_{cfl}$ &  $\overline{\Omega}$ & $\Omega$  & $\overline{N_0}$  & $ N_0$   \\ \hline
    0.7 & 0.2  & $[-400,0]$  & $[0, 400]$  & 2 & 2 
  \end{tabular}
  \caption{Parameters for test case 1.}\label{tab:params-case-1}
\end{table}
The $C_{cfl}$ constant is set such that the timestep is  five times less than the theoretically stable bound in order to predict the variation of the characteristic speed coming (solely) from the coupling interface. 

\begin{figure}[!htb]
  \centering
  \includegraphics[width=0.32\textwidth]{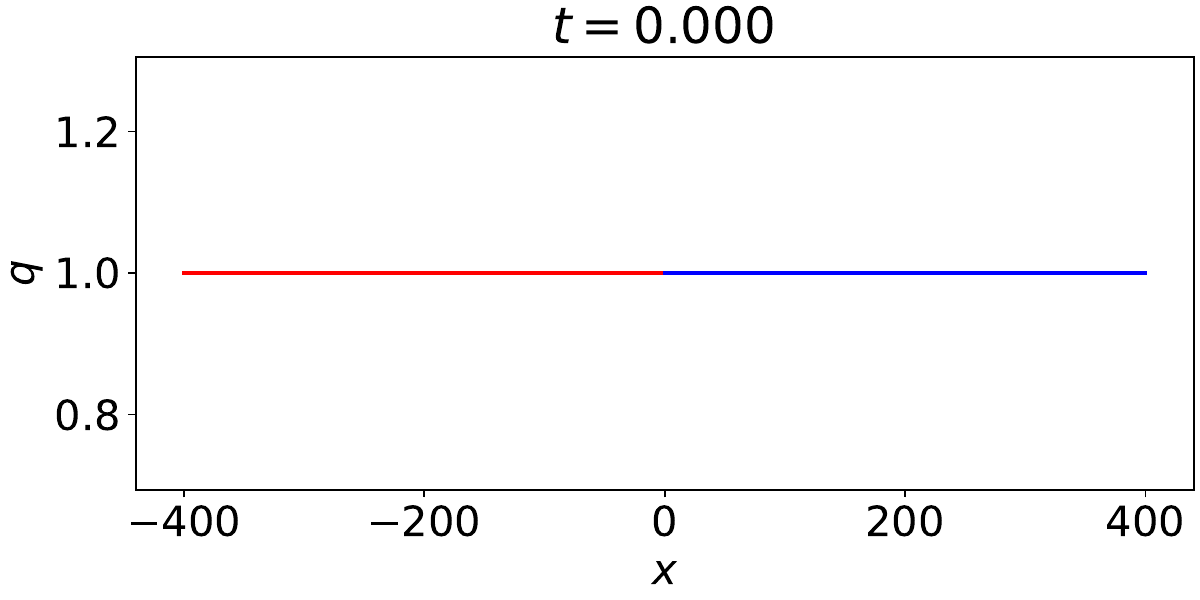}
  \includegraphics[width=0.32\textwidth]{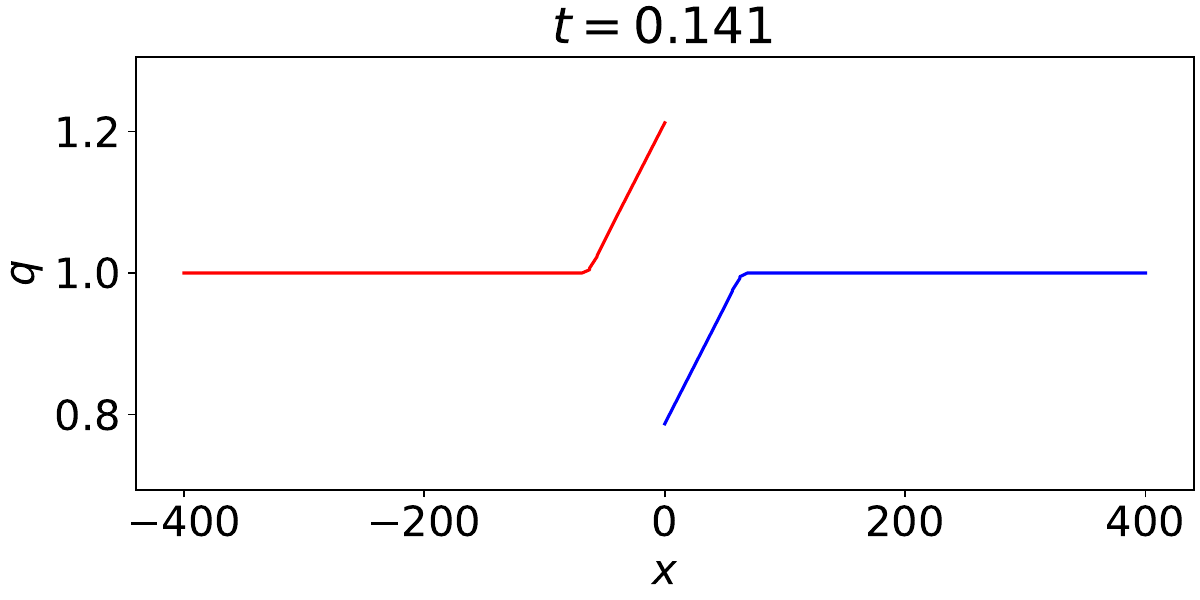}
  \includegraphics[width=0.32\textwidth]{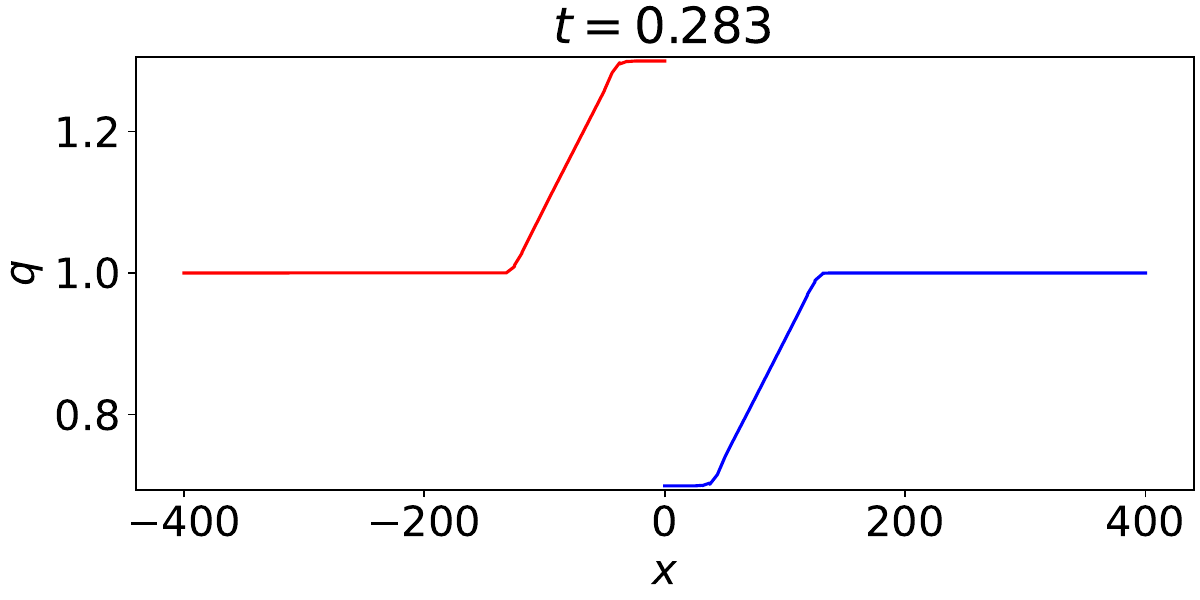}
  \includegraphics[width=0.32\textwidth]{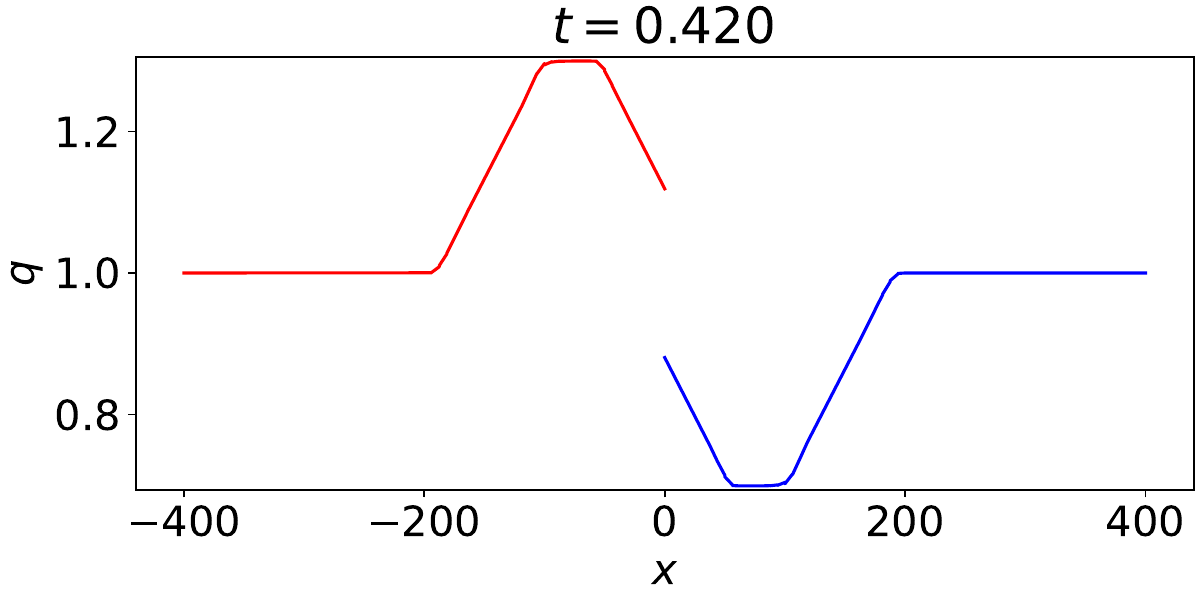}
  \includegraphics[width=0.32\textwidth]{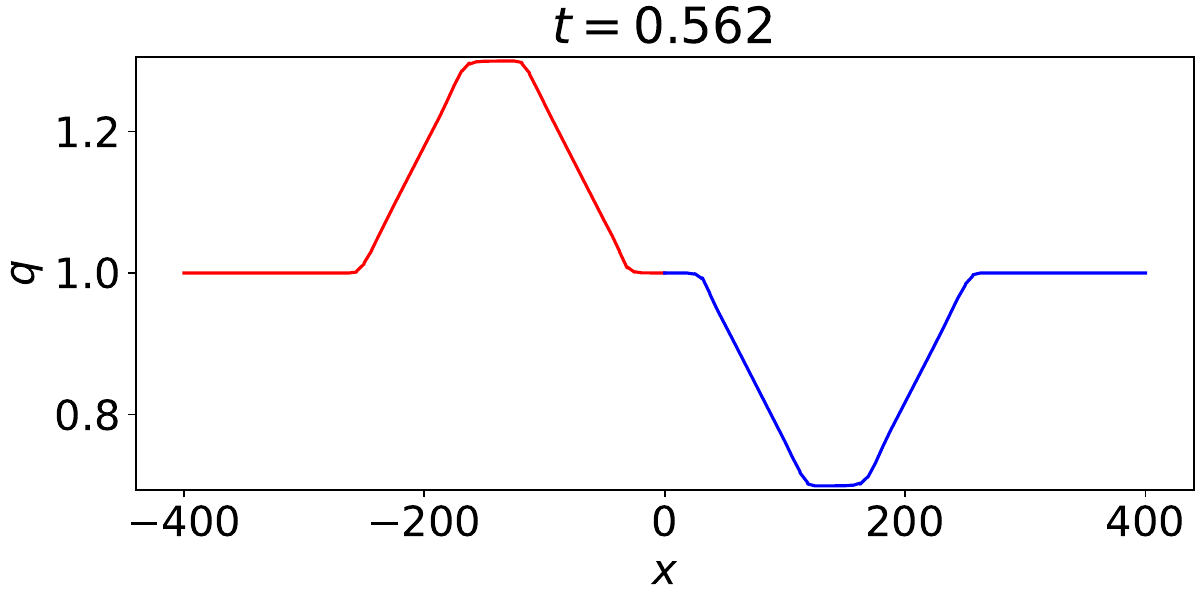}
  \includegraphics[width=0.32\textwidth]{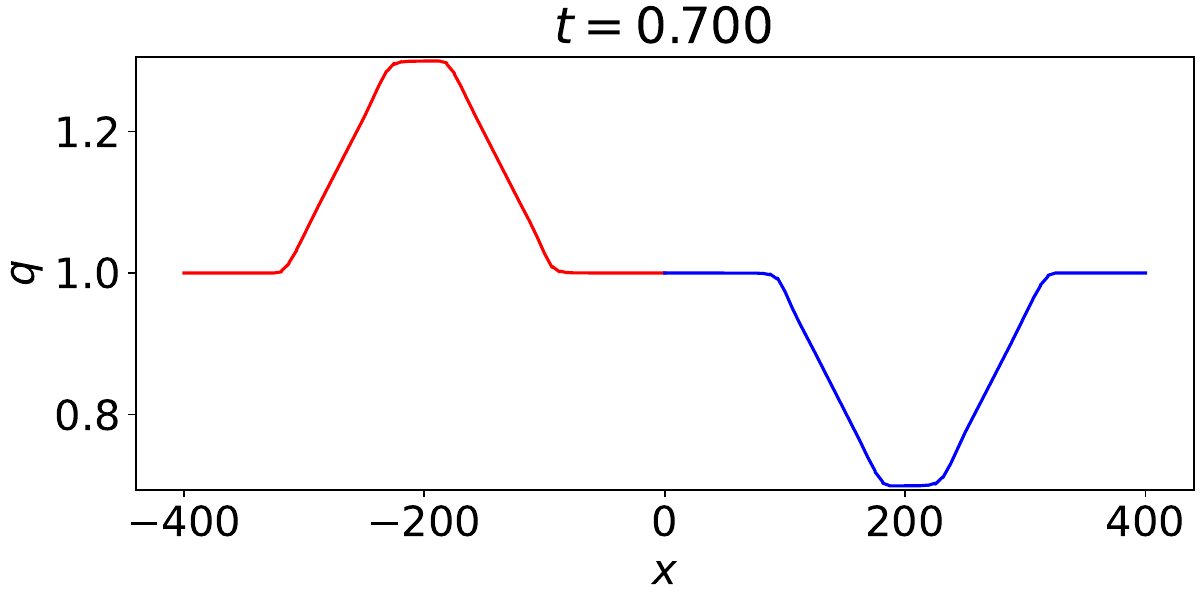}
  \caption{Time series of the momentum $q$ for Case 1 on level $L=5$.}
  \label{fig:lin-outtake-time-series}
\end{figure}

The behaviour of the resulting momentum is illustrated in Figure~\ref{fig:lin-outtake-time-series} and the error convergence rates are listed in the second and third column of Table~\ref{tab:errors}. The results match those obtained by the classical RP-coupling strategy in~\cite{hertyCouplingCompressibleEuler2019}, however, due to  appearance of kinks the convergence rates are expectably around 1.6. 

\begin{figure}[!htb]
  \centering
  \includegraphics[width=0.32\textwidth]{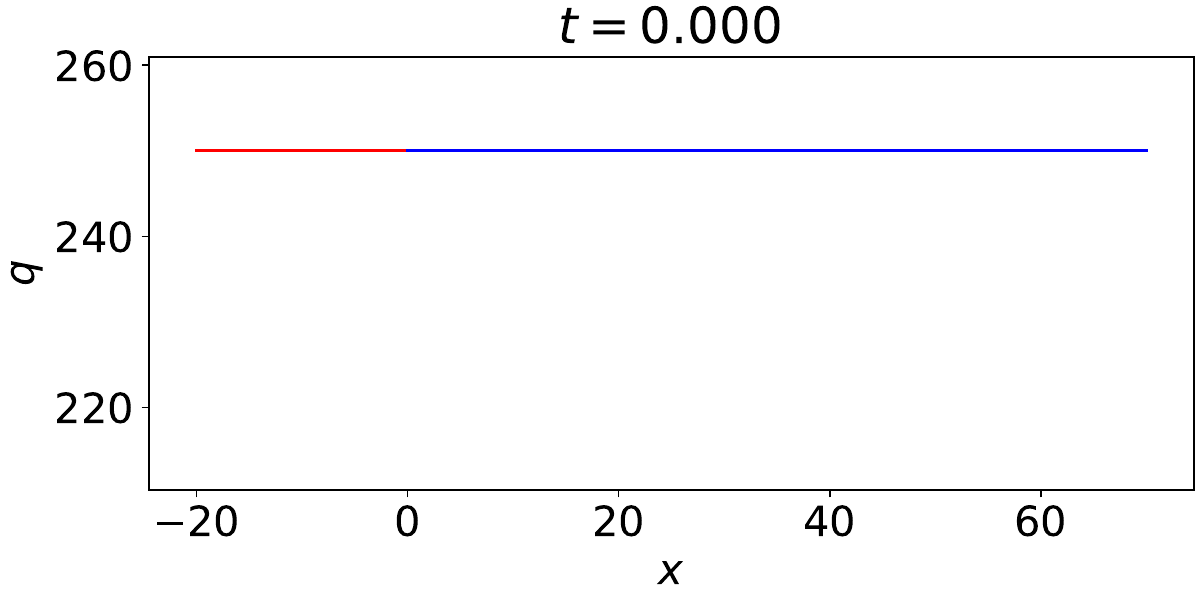}
  \includegraphics[width=0.32\textwidth]{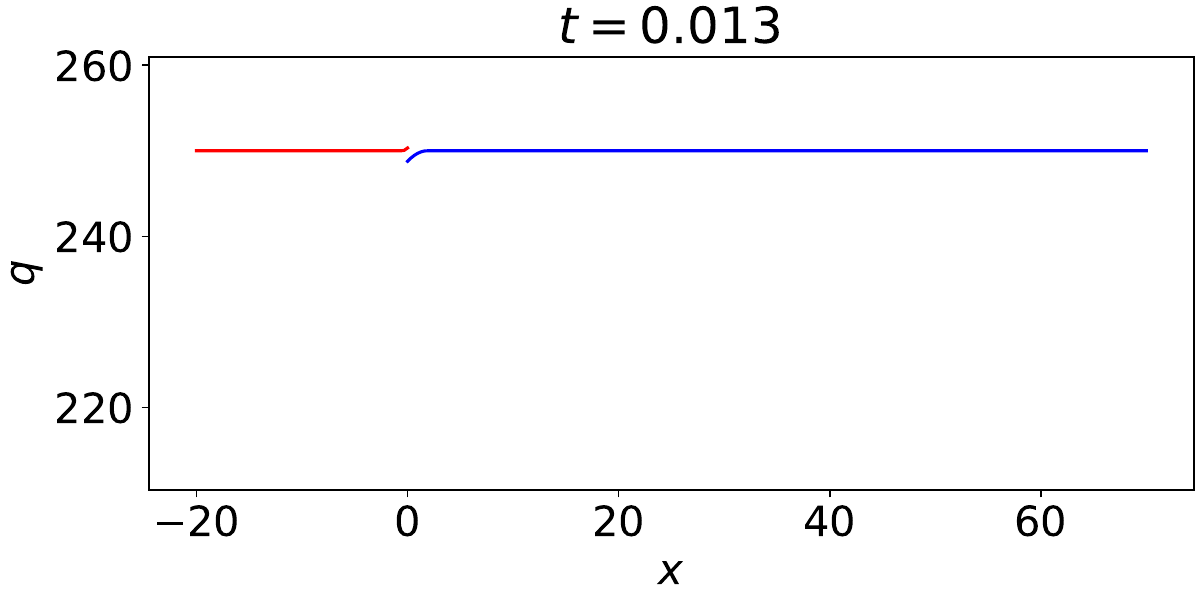}
  \includegraphics[width=0.32\textwidth]{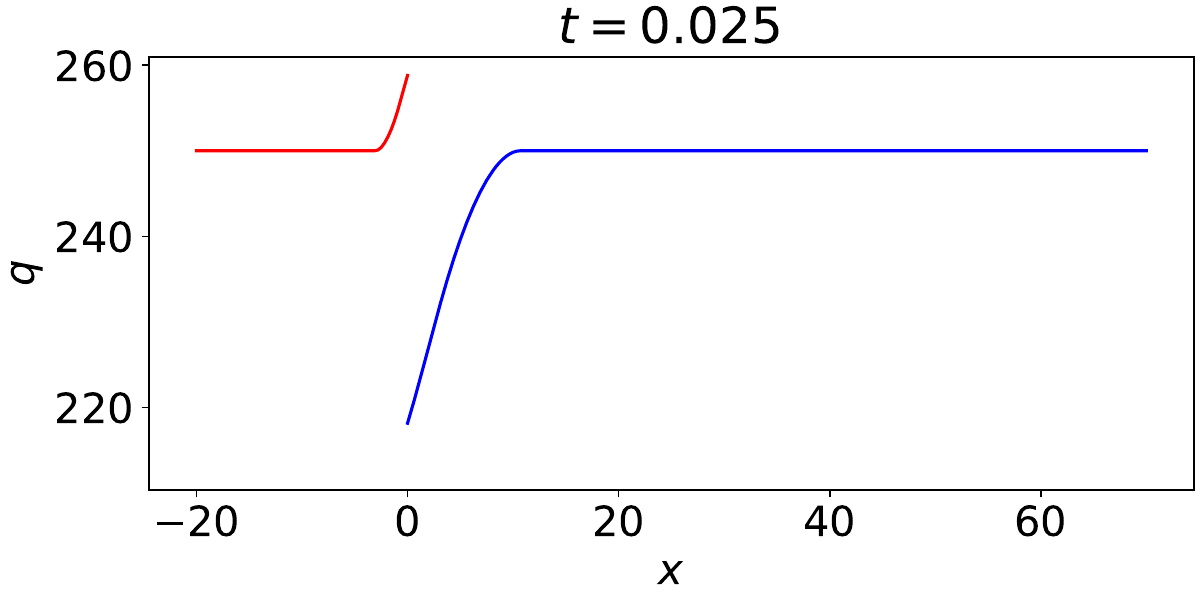}
  \includegraphics[width=0.32\textwidth]{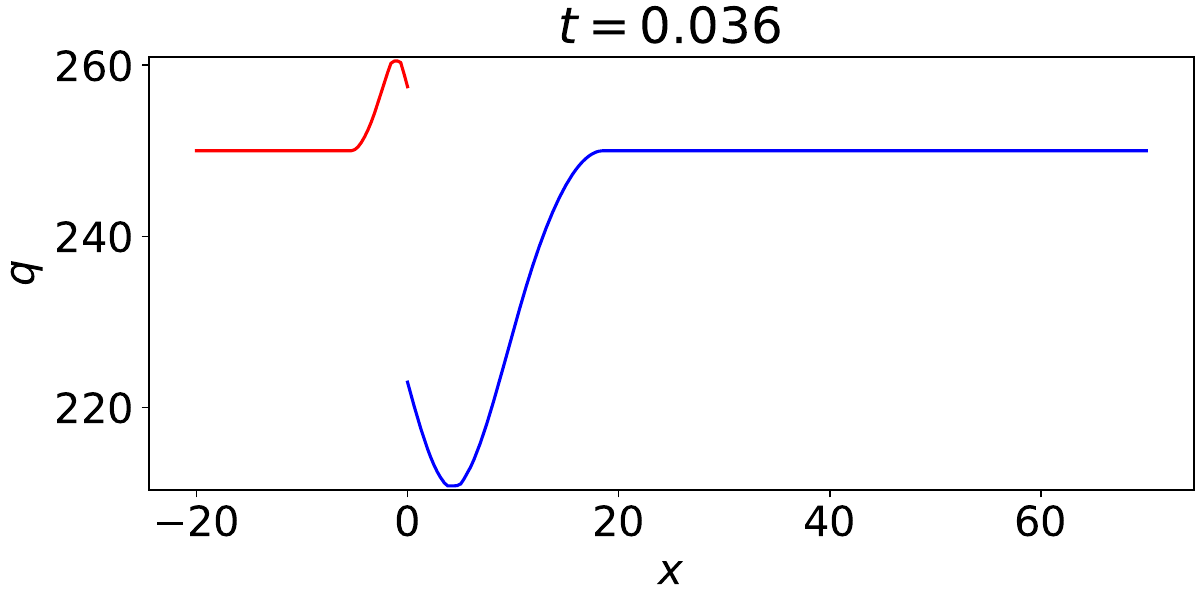}
  \includegraphics[width=0.32\textwidth]{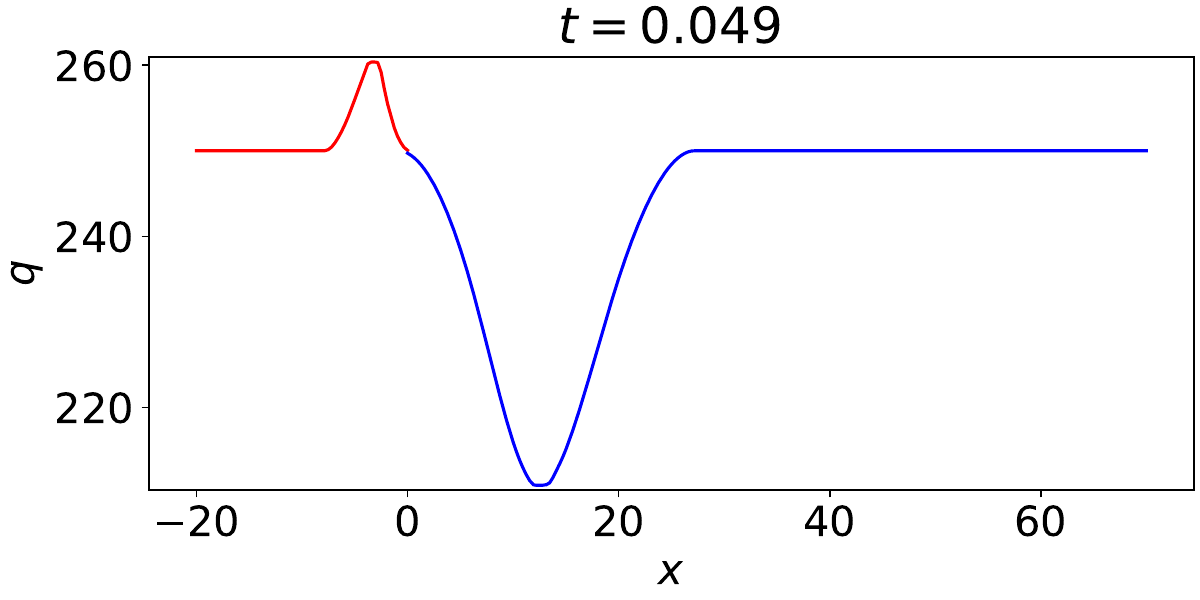}
  \includegraphics[width=0.32\textwidth]{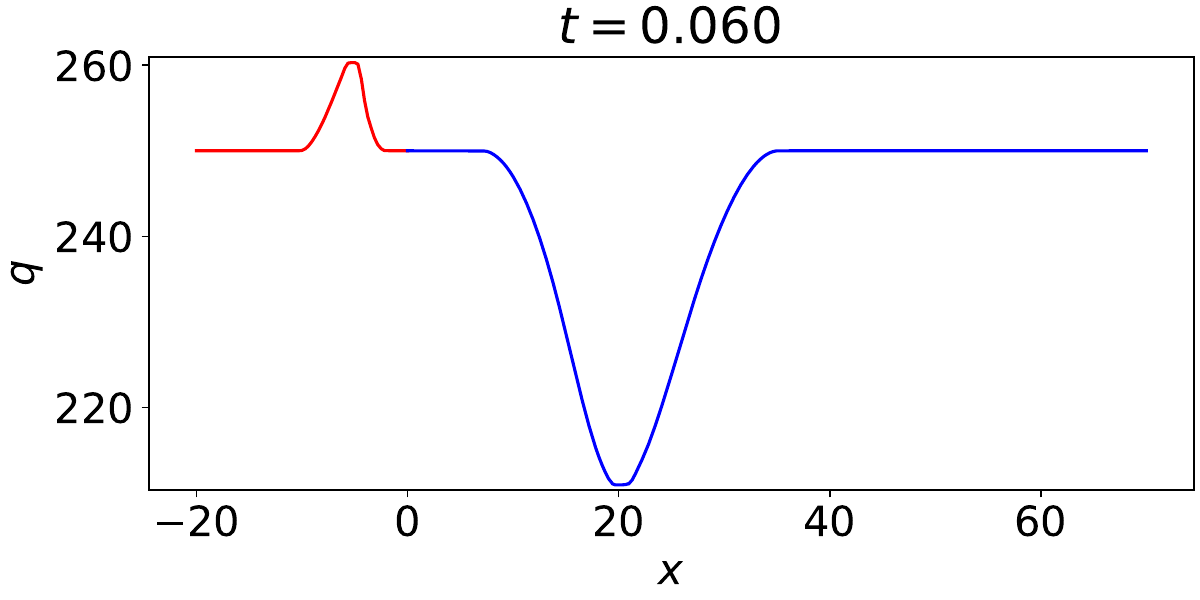}
  \caption{Time series of the momentum $q$ for Case 2 on level $L=5$.}
  \label{fig:nonlin-outtake}
\end{figure}

\paragraph{Case 2: Nonlinear Outtake}
We set the initial condition in both pipes to  a constant state $\rho(0,x) = 1.0$, $v(0,x) = 250.0$ and $p(0,x) = 146820.4$. For the outtake we prescribe a smooth perturbation given by a periodic cubic spline through the points listed in Table~\ref{tab:case2-spline}.
\begin{table}[!htb]
  \centering
\begin{tabular}{l|lllll}
$t$ &  0.01& 0.02& 0.03& 0.04& 0.05  \\\hline
$\mathcal{E}(t)$  &  0.0 & 20.0 & 50.0 & 20.0 & 0.0  
\end{tabular}
\caption{Spline nodes and values for test case 2.}\label{tab:case2-spline}
\end{table}
The setup of the domain and its discretization parameters are summarized in Table~\ref{tab:params-case-2} and the visualization is given in Figure~\ref{fig:nonlin-outtake}.
\begin{table}[!htb] 
  \centering 
  \begin{tabular}{c|c|cc|cc}
    $T$ & $C_{cfl}$ &  $\overline{\Omega}$ & $\Omega$  & $\overline{N_0}$  & $ N_0$   \\ \hline
    0.06 & 0.2  & $[-20,0]$  & $[0, 70]$  & 2 & 7 
  \end{tabular}
  \caption{Parameters for test case 2.}\label{tab:params-case-2}
\end{table}
Again the timestep is reduced in order to account for varaitions of the characteristic speed originating from the coupling  boundary.

Table~\ref{tab:errors} shows EoC rates of about two for this test case. This is a reasnable rate since the entire variation of the solution enters the domain from the coupling interface that, in turn, is approximated with second order in time.

\paragraph{Case 3: Density bump passing a nonlinear outtake}
In the final test case we consider a perturbation in the density of the constant initial state  $\rho(0,x) = 1.0 + s(x)$, $v(0,x) = 30.0$ and $p(0,x) = 146820.4$ where $s$ is a cubic periodic spline on the points detailed in Table~\ref{tab:case3-spline}.
\begin{table}[!htb]
  \centering
\begin{tabular}{l|rrrrr}
$x$ & $-15$& $-12.5$& $-10$ & $-7.5$& $-5$  \\\hline
$s(x)$  &  0.0& 0.1& 0.2& 0.1& 0.0
\end{tabular}
\caption{Spline nodes and values for test case 2.}\label{tab:case3-spline}
\end{table}
The outtake is set to $\mathcal{E}(t) = 3$ and the domain as well as the mesh are specified in Table~\ref{tab:params-case-3}.
\begin{table}[!htb] 
  \centering 
  \begin{tabular}{c|c|cc|cc}
    $T$ & $C_{cfl}$ &  $\overline{\Omega}$ & $\Omega$  & $\overline{N_0}$  & $ N_0$   \\ \hline
    0.6 & 0.9  & $[-20,0]$  & $[0, 20]$  & 2 & 2 
  \end{tabular}
  \caption{Parameters for test case 3.}\label{tab:params-case-3}
\end{table}
Note that, in contrast to the previous ones, in this case the timestep size is dominated by the solution in the interior of the domains and, thus, $C_{cfl}= 0.9$ is sufficient for a stable discretization.

\begin{figure}[!htb]
  \centering
  \includegraphics[width=0.32\textwidth]{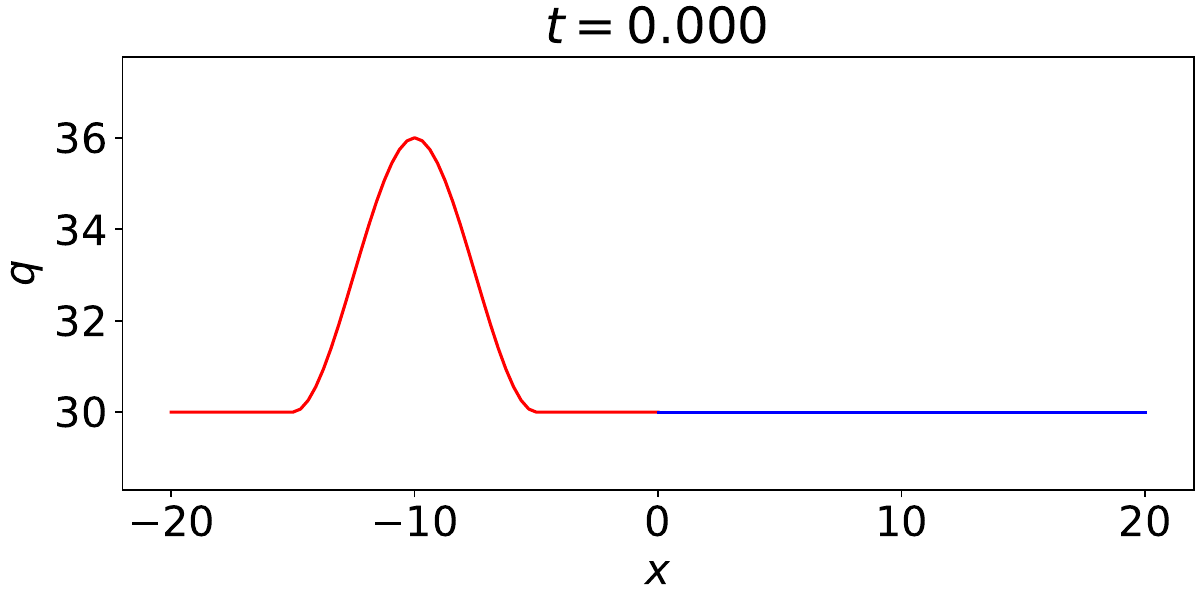}
  \includegraphics[width=0.32\textwidth]{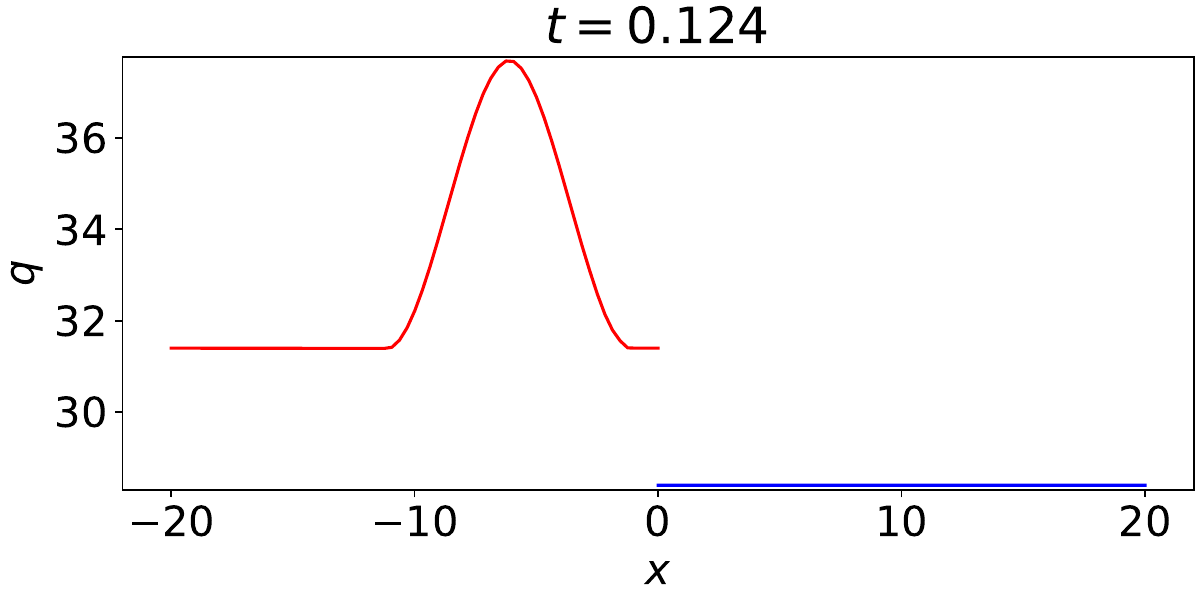}
  \includegraphics[width=0.32\textwidth]{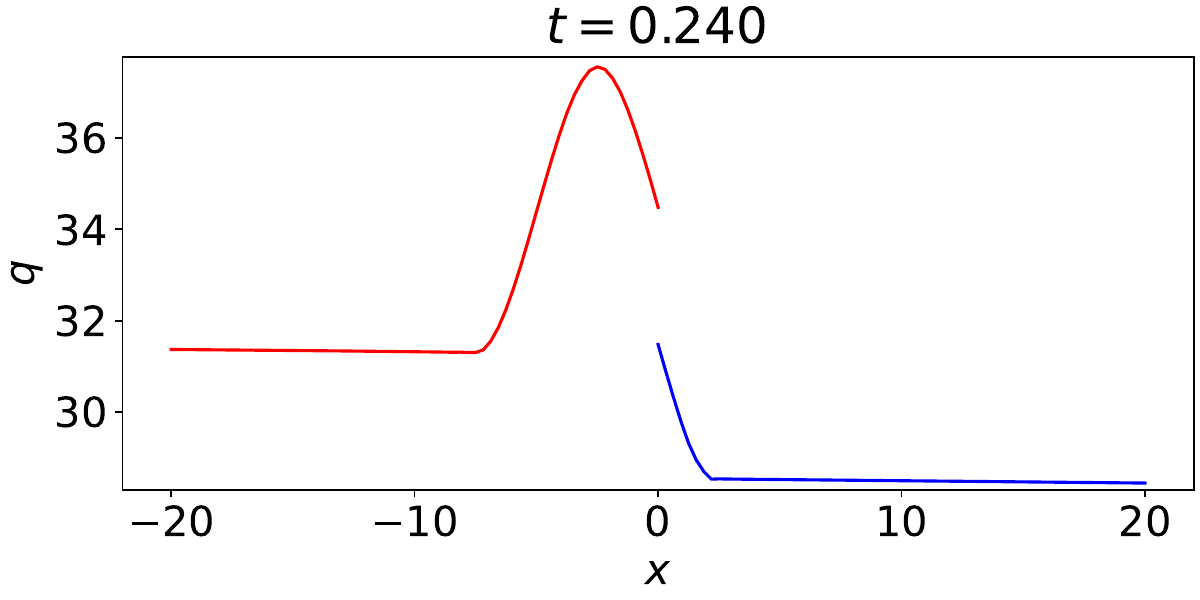}
  \includegraphics[width=0.32\textwidth]{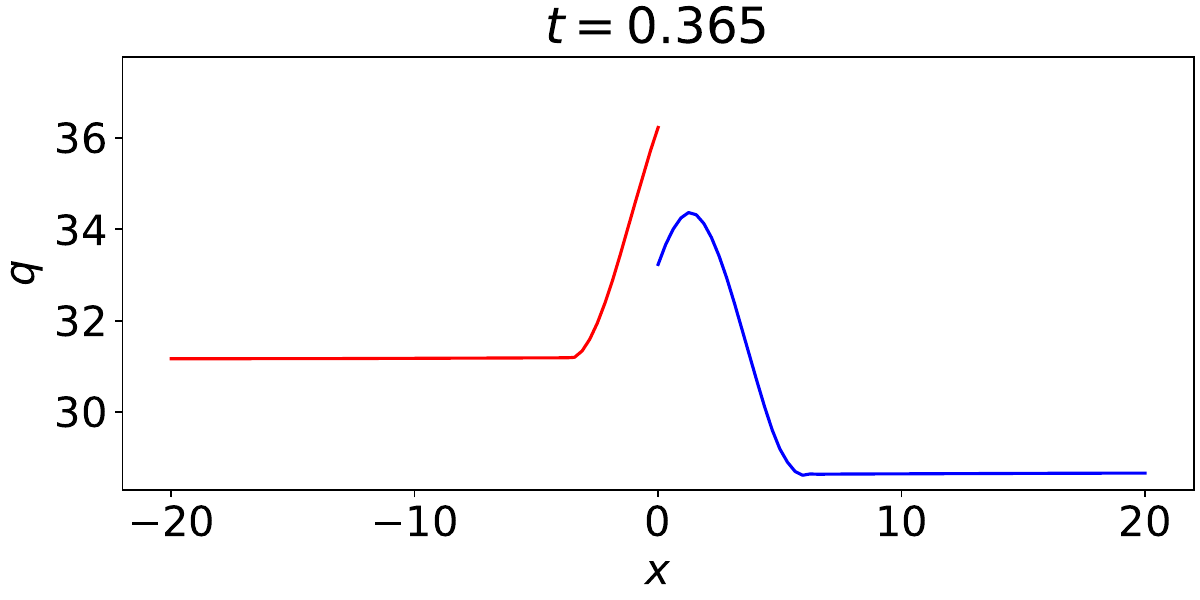}
  \includegraphics[width=0.32\textwidth]{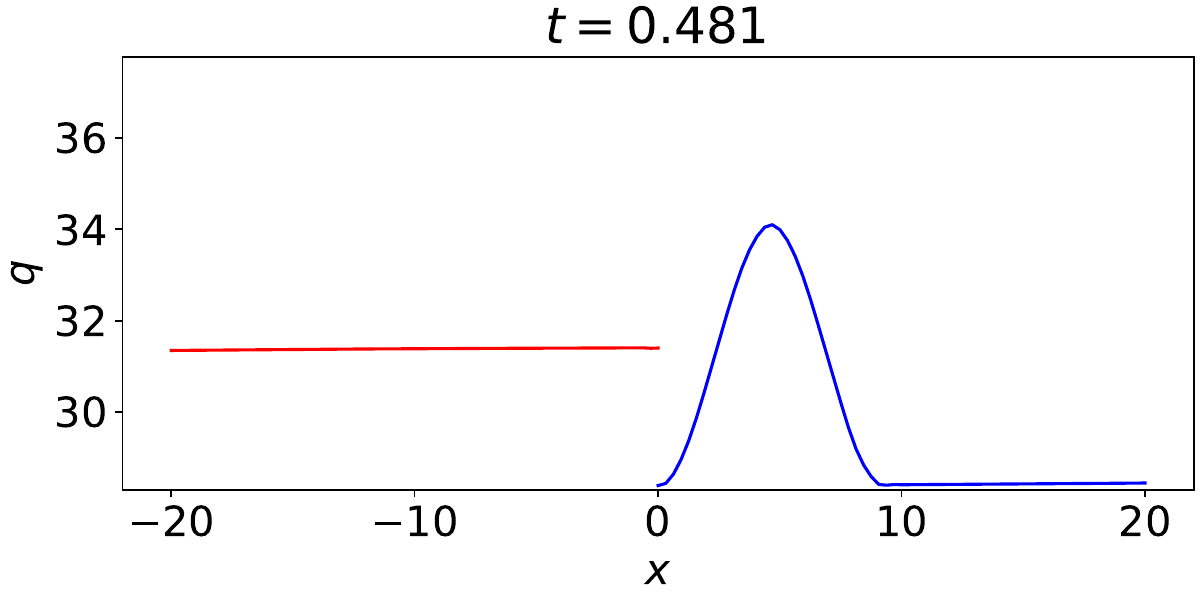}
  \includegraphics[width=0.32\textwidth]{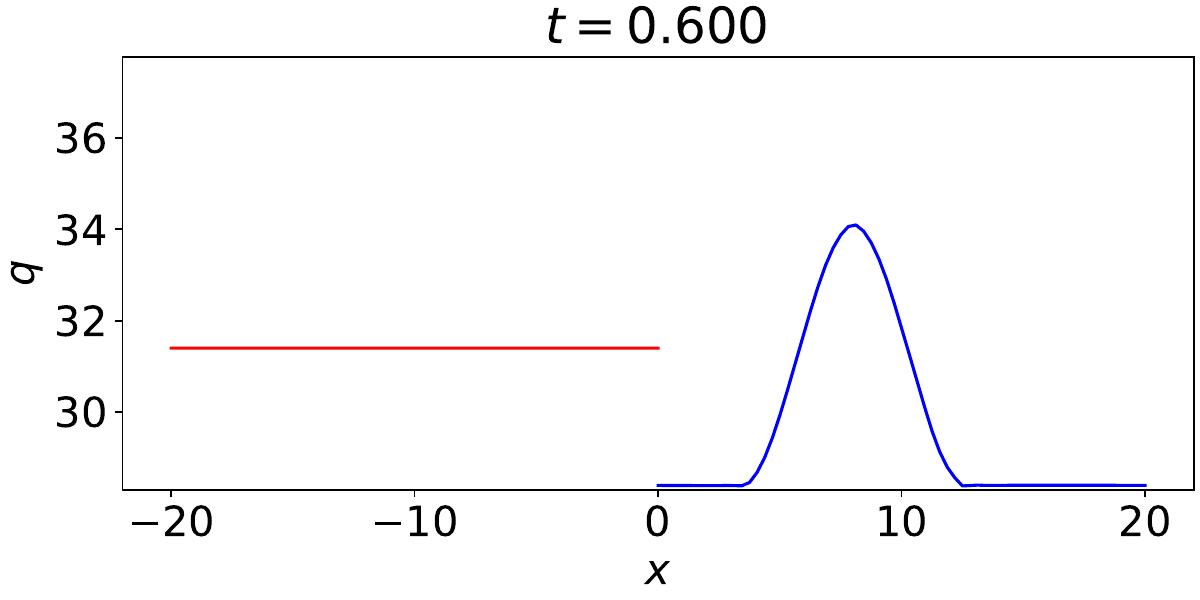}
  \caption{Time series of the momentum $q$ for Case 3 on level $L=5$.}
  \label{fig:perturbation-across-jump}
\end{figure}

The instant jump in the outtake at $t=0$ leads to three outgoing waves that quickly leave the respective coupled domains. Since the velocity and pressure are constant in each of the domains, the density perturbation $s$ is transported across the coupling interface (which contains a jump in momentum). This behaviour is shown in Figure~\ref{fig:perturbation-across-jump}. The convergence rates listed in Table~\ref{tab:errors} are around 2.7 which is surprisingly high. This might be explained by the dominance in the errors of the solution in domain interiors, in contrast to Case 2 where the error was clearly dominated by the coupling boundary.   

\begin{table}[]
  \centering
\begin{tabular}{l|ll|ll|ll}
    & \multicolumn{2}{l|}{Case 1} & \multicolumn{2}{l|}{Case 2} & \multicolumn{2}{l}{Case 3} \\
$L$ &  $err$     & EoC      & $err$     & EoC      & $err$     & EoC      \\ \hline
3   & 7273.33       &          & 29229.9         &          & 51.4209         &          \\
4   & 2112.71       & 1.78     & 8450.18         & 1.79     & 8.40335         & 2.61     \\
5   & 626.409       & 1.75     & 1840.29         & 2.20     & 1.30399         & 2.69     \\
6   & 211.697       & 1.57     & 393.232         & 2.23     & 0.213337        & 2.61     \\
7   & 69.3842       & 1.61     & 97.0706         & 2.02     & 0.0238593       & 3.16    
\end{tabular}
\caption{Errors and estimated order of convergence for all test cases.}\label{tab:errors}
\end{table}

\begin{rem}
  Algorithm~\ref{alg:CGRP} yields a solution to the coupling problem that is  of second order in time and might reduce the overall order of a higher resolution scheme to two, as we saw in Case 2. One way to overcome  this problem is deriving a higher-order coupled GRP. Another possibility is adding a sufficient number of time-synchronisation points in the interior of each timestep and solving multiple coupled GRPs in between.  In conjunction with dense output techniques, see e.g.~\cite{ketcheson2017dense}, this yields a high-order coupling method with essentially independent explicit time-discretizations.
\end{rem}

\section*{Acknowledgments}
Aleksey Sikstel acknowledges funding from the German Science Foundation DFG through the research unit “SNuBIC” (DFG-FOR5409).
